\documentclass[12pt]{article}

\usepackage[T2A]{fontenc}
\usepackage[cp866]{inputenc}

\usepackage{amsfonts}
\usepackage{amssymb}
\usepackage{amsmath}
\usepackage{amsthm}

\def \le {\leqslant}
\def \ge {\geqslant}

\begin{document}

\begin{Large}
\centerline{\bf On small fractional parts of polynomials}
\end{Large}
 \vskip+1.5cm \centerline{\bf Moshchevitin N.G. \footnote{ Research is supported by
grants RFFI 06-01-00518, MD-3003.2006.1, NSh-1312.2006.1 and INTAS 03-51-5070 }} \vskip+1.5cm

{\bf Abstract.}\,\,
 We prove that  for any real polynomial $f(x) \in\mathbb{R} [x]$  the set
$$
\{ \alpha \in \mathbb{R}:\,\,\, \liminf_{n\to \infty}\,\, n\log n \, ||\alpha f(n)||\,
>0\}
$$
 has positive Hausdorff dimension.
 Here $||\xi ||$ means the distance from $\xi $ to the nearest integer.
Our result is based on an original method due to Y. Peres and W. Schlag.
 \vskip+1.5cm

{\bf 1.\,\, Introduction.} In this paper we prove the following result.

{\bf Theorem 1.}\,\,{\it Let $t_n$ be an increasing sequence of reals, $ t_0 \ge 1$. Let
\begin{equation}
\frac{t_{n+1}}{t_n} = 1+\frac{\gamma }{n} +  O\left( \frac{1}{n^{1+\varepsilon_1}}\right),\,\,\, \gamma > 0, \varepsilon_1 > 0, \,\,\, n\to
+\infty . \label{LACC}
\end{equation}
Then the Hausdorff dimension of the set
$$
\{ \alpha \in \mathbb{R}:\,\,\, \liminf_{n\to \infty}\,\, n\log n \, ||\alpha t_n||\,
>0\}
$$
is greater than $\frac{\gamma}{\gamma+1}$.}

As for a polynomial $ f(x)\in\mathbb{R} [x]$  of degree ${\rm deg } f = d\ge 1$ one has
$$
\frac{f(n+1)}{f(n)} = 1+\frac{d }{n} +  O\left( \frac{1}{n^{2}}\right),\,\,\,n\to +\infty
$$
we immediately obtain the following

{\bf Corollary.}\,\,{\it Let $f(n)$ be a polynomial of degree $d\ge 1$.  Then the Hausdorff dimension of the set
$$
\{ \alpha \in \mathbb{R}:\,\,\, \liminf_{n\to \infty}\,\, n\log n \, ||\alpha f(n)||\,
>0\}
$$
is greater than $\frac{d}{d+1}$.}

We would like to note that for $d=1$ our result is trivial. Consider the case $ d=2$. In \cite{S} W. Schmidt asks if it is true that for any
real $\alpha$ we have
$$
\liminf_{n\to +\infty} n\log n \, ||\alpha n^2|| =0.
$$
Our Theorem 1 gives a negative answer to this question. The complete proof of  Theorem  1  is given in Sections 2,3,4 below. We would like to
say that our result is based on the original construction introduced recently by Y. Peres and W. Schlag \cite{P}.

Consider a polynomial $f(n) = a n^k$ with a nonzero integer coefficient $a$  and a function $\psi (n)$ decreasing to zero as $ n\to +\infty$.
From Cassels' theorem \cite{C} it follows that in the case when
$$
\sum_{n=1}^\infty \psi (n) = +\infty
$$
for almost all real $\alpha $ we have
$$
\liminf_{n\to +\infty} \,\, (\psi (n))^{-1}\, ||\alpha f(n)|| =0.
$$
So for almost all real numbers $\alpha$ we have
$$
\liminf_{n\to +\infty}  n\log n\, ||\alpha f(n)|| =0
$$
and
 the set constructed in  Theorem 1 has Lebesgue measure zero. Also we would like to note that in
the case
$$
\sum_{n=1}^\infty \psi (n) < +\infty
$$
for almost all real $\alpha $ we have
$$
\liminf_{n\to +\infty}  \,\, (\psi (n))^{-1}\, ||\alpha f(n)|| >0.
$$

 There is a lot of results on upper bounds for small
fractional parts of polynomials (see \cite{B}).
 For the polynomial $ f(n) = n^2$ we would like to mention here the result due to A. Zaharescu \cite{Z}:
 for $\varepsilon > 0$ and every  real $\alpha $
there are infinitely many natural numbers $n$ such that
$$
||\alpha n^2|| \le c(\varepsilon ) n^{-\frac{2}{3}+\varepsilon},
$$
where $c(\varepsilon) $ is a positive constant.

 {\bf 2.\,\, Parameters. }

One can easily see that under the condition (\ref{LACC})  for $ m > n$ we have the following estimate:
\begin{equation}
\frac{t_{m}}{t_n} = \prod_{j = n}^m \left( 1+\frac{\gamma}{n} +  O\left( \frac{1}{n^{1+\varepsilon_1}}\right)\right) = \exp \left(\sum_{j = n}^m
\frac{\gamma}{n} +  o(1) \right) \asymp \left( \frac{m}{n}\right)^\gamma
 .
 \label{EST}
\end{equation}
When $ n = 1$ we get
\begin{equation}
t_m \asymp m^\gamma
 .\label{T}
\end{equation}
For $ n > 1$ define
\begin{equation}
h (n) = n^{1+\frac{1}{\gamma}}\log^{\frac{2}{\gamma}}n,\,\,\, \delta_n  = \frac{1}{cn\log n},\,\,\, c = 60 \log \left( 2 +\frac{1}{\gamma }
\right).\label{PARA}
\end{equation}
From (\ref{EST}) in follows that for   $ m \ge h(n)$ one has
$$
\frac{t_{m}}{t_n}\gg n\log^2 n
$$
and when $n$ is large enough one has
\begin{equation}
\frac{t_{m}}{t_n} \ge \frac{1}{\delta_n}
 .\label{D}
\end{equation}
Let $ n_0 \ge 2$. Define natural numbers $ n_j$ inductively by the equality
\begin{equation}
n_{j+1} = h ({n_j}). \label{NJ}
\end{equation}
For a given $\varepsilon_2 > 0$ we can choose  $ n_0$ to be  large enough so that
\begin{equation}
n_{k-1}^{1+\frac{1}{\gamma}} \le n_k \le n_{k-1}^{1+\frac{1}{\gamma}+\varepsilon_2},\,\,\, \forall k \in \mathbb{N} \label{nei}
\end{equation}
and
\begin{equation}
n_{0}^{\left(1+\frac{1}{\gamma}\right)^k} \le n_k \le n_{0}^{\left(1+\frac{1}{\gamma}+\varepsilon_2\right)^k},\,\,\, \forall k \in \mathbb{N}.
\label{n}
\end{equation}
As a corollary we immediately have
\begin{equation}
\log n_k \asymp \log n_{k-1}
 .\label{log}
\end{equation}
For $ v \in (0,1)$ we consider the series
\begin{equation}
\sum_{k= 2}^\infty 3^k \, \frac{t_{n_k}^v       }{ t_{n_{k-1}} }\, \frac{\delta_{n_{k-1}}}{\delta_{n_k}^v}.
 \label{ser}
\end{equation}
From(\ref{T}, \ref{nei}, \ref{log}) it follows that
$$
\frac{t_{n_k}^v       }{ t_{n_{k-1}} }\, \frac{\delta_{n_{k-1}}}{\delta_{n_k}^v} \ll n_{k-1}^\omega
$$
where
$$
\omega = \left(\left( 1+\frac{1}{\gamma} +\varepsilon_2\right) v - 1\right)(\gamma +1).
$$
We should note that in the case $ v \le \frac{\gamma}{\gamma +1}$ one can  take small positive $\varepsilon_2$ to make $ \omega <0$. Then from
the lower bound in (\ref{n}) it follows that the series (\ref{ser}) converges.

{\bf 2.\,\, Sets. }

Define
\begin{equation}
E_n = \bigcup_{a=0}^{\lceil t_n\rceil} \left[ \frac{a}{t_n} - \frac{\delta_n}{t_n}, \frac{a}{t_n} + \frac{\delta_n}{t_n} \right]\bigcap [0,1].
\label{E}
\end{equation}
Let
\begin{equation}
l_n = \left\lfloor \log_2 \frac{2t_n}{\delta_n} \right\rfloor,\,\,\, 2^{l_n}\le \frac{2t_n}{\delta_n} < 2^{l_n+1}. \label{L} \end{equation}
 Each
segment form the union (\ref{E}) can be covered by a dyadic interval of the form
$$
\left( \frac{b}{2^{l_n}}, \frac{b+z}{2^{l_n}}\right),\,\,\, z = 1,2 .$$ Let $A_n$ be the smallest union of all such dyadic segments which cover
the whole set $E_n$. Define
$$
A_n^c = [0,1] \setminus A_n = \bigcup_{\nu = 1}^{\tau_\nu } J_\nu
$$
where closed  segments $J_\nu $ are of the form
\begin{equation}
\left[ \frac{b}{2^{l_n}}, \frac{b+1}{2^{l_n}}\right] . \label{J}
\end{equation}
Define
$$
B_n = \bigcap_{j\le n} A_j^c,\,\,\, B = \bigcap_{j=1}^\infty A_j^c.
$$
Each $B_n$ can be written in the form $B_n= \bigcup_{\nu = 1}^{T_\nu } J_\nu^{(n)}$ where $ J_\nu^{(n)}$ are of the form (\ref{J}).

For our purpose it is enough to prove that the Hausdorff dimension of the set $B$ is  not less than $  \frac{\gamma}{\gamma +1}. $

 {\bf 3.\,\, Lemmata. }

 {\bf  Lemma 1. }
\,\, {\it Let $ B_n = \bigcup_{\nu = 1}^{T_n} J_\nu ^{(n)} \neq \varnothing$ and $ m \ge h(n)$. Then for every $\nu$ from the interval $ 1\le
\nu \le T_n$ one has
$$
\mu \left( J_\nu^{(n)} \bigcap A_m \right) \le 5 \delta_m \mu \left( J_\nu^{(n)} \right).
$$}

Proof. Write $A_m = \cup_i I_i$ where each $ I_i$ is of the form $ (b/2^{l_m}, (b+1)/2^{l_m}).$ Let $I_i \cap J_\nu ^{(n)} \neq \varnothing$ and
$J_\nu ^{(n)} = [a/2^{l_n}, (a+1)/2^{l_n}].$ Then for some natural $z$ we have
$$
\frac{z}{t_m} \in \left[ \frac{a}{2^{l_n}} - \frac{1}{2^{l_m}},
 \frac{a+1}{2^{l_n}} + \frac{1}{2^{l_m}}\right].
 $$
 For each $\nu$ the quantity $z$ can take not more than
 $$
 W = \left\lfloor  \left( \frac{1}{2^{l_n} } + \frac{2}{2^{l_m}} \right) t_m \right\rfloor
+ 1
$$
values. Now
\begin{equation}
\mu \left( J_\nu^{(n)} \bigcap A_m \right) \le \mu\left( I_i\right) W = \frac{2^{l_n}}{2^{l_m}} W \mu  \left( J_\nu ^{(n)}\right) \le \mu
\left(J_\nu ^{(n)}\right) \left( 3\,\frac{t_m}{2^{l_m}} + \frac{2^{l_n}}{2^{l_m}}\right). \label{1}
\end{equation}
From the definition (\ref{L}) of $l_n$ it follows that
\begin{equation}
\frac{t_m}{2^{l_m}} \le \delta_m
 . \label{2}
\end{equation}
From  (\ref{L}) and (\ref{D}) it follows that
\begin{equation}
\frac{2^{l_n}}{2^{l_m}}\le  2\delta_m.
 . \label{3}
\end{equation}
Substituting (\ref{2},\ref{3}) into (\ref{1}) we get  Lemma 1.

{\bf  Lemma 2. } \,\, {\it Let $m = h(n), M = h(m)$ and let $n$ be large enough. Consider a segment $J_\nu^{(n)} =
\left[\frac{a}{2^{l_n}},\frac{a+1}{2^{l_n}}\right]$ from the partition
\begin{equation}
\varnothing \neq B_n =\bigcup_{\nu= 1}^{T_n} J_\nu^{(n)} . \label{union} \end{equation}
 Suppose that we know that
\begin{equation}
\mu \left( J_\nu^{(n)} \bigcap B_m \right) \ge \frac{1}{2} \,\, \mu \left( J_\nu^{(n)} \right)
 \label{ini}
\end{equation}
(in fact it means that $B_m \neq \varnothing$). Then
$$
\mu \left( J_\nu^{(n)} \bigcap B_M \right) \ge \frac{5}{6} \,\, \mu \left( J_\nu^{(n)}\bigcap B_m \right) .$$ }

Proof.\,\, We have
$$
 J_\nu^{(n)}\bigcap B_M
=\left( ...\left(\left(
 J_\nu^{(n)}\bigcap B_m\right) \setminus A_{m+1}\right) \setminus  ...\right) \setminus A_M
 $$
 and
 \begin{equation}
\mu \left( J_\nu^{(n)} \bigcap B_M \right) = \mu \left( J_\nu^{(n)} \bigcap B_m \right) - \sum_{j = m+1}^M \mu \left( \left( J_\nu^{(n)} \bigcap
B_m \right)\bigcap A_j \right).
 \label{1w}
 \end{equation}
Define
$$
n(j) = \max \{ k:\,\, h(k) \le j\}.
$$
As m = h(n) < j we see that $ n(j) \ge n$ and
\begin{equation}
B_{n(j)} \subseteq B_n. \label{Aw} \end{equation}
On
 another hand  $ h(m) = M \ge j $ and for any $m'\ge m$ we have $h(m') \ge M \ge j$. So $ n(j)
\le m$ and
\begin{equation}
B_m \subseteq B_{n(j)} . \label{Bw} \end{equation} From (\ref{Bw}) it follows that
$$
 J_\nu^{(n)}\bigcap B_m
 \subseteq  J_\nu^{(n)}\bigcap B_{n(j)} =
 \bigcup_{\nu_1:\,\, J_{\nu_1}^{(n(j))} \subseteq  J_\nu^{(n)}} J_{\nu_1}^{(n(j))}.
 $$
 So
$$
\left( J_\nu^{(n)}\bigcap B_m \right) \bigcap A_j \subseteq \left( J_\nu^{(n)}\bigcap B_{n(j)} \right) \bigcap A_j
=
 \bigcup_{\nu_1:\,\, J_{\nu_1}^{(n(j))} \subseteq  J_\nu^{(n)}} \left(J_{\nu_1}^{(n(j))}\bigcap A_j \right).
 $$
 As $ j \ge h(n(j))$ we can apply Lemma 1. From Lemma 1 it follows that
 $$
 \mu \left( J_{\nu_1}^{(n(j))} \bigcap A_j
 \right)\le 5 \delta_j \mu \left( J_\nu^{(n(j))} \right).
 $$ Summation over $\nu_1$ gives
 \begin{equation}
\mu \left(\left( J_\nu^{(n)}\bigcap B_{n(j)} \right) \bigcap A_j\right) \le 5 \delta_j \mu \left( J_\nu^{(n)} \bigcap B_{n(j)}\right).
 \label{2w}
 \end{equation}
Moreover from (\ref{Aw}) it follows that
 \begin{equation}
\mu \left( J_\nu^{(n)}\bigcap B_{n(j)} \right) \le \mu \left( J_\nu^{(n)} \right).
 \label{3w}
 \end{equation}
 From (\ref{ini}) we have
 \begin{equation}
  \mu \left( J_\nu^{(n)} \right) \le 2
 \mu \left( J_\nu^{(n)} \bigcap B_m \right).
 \label{4w}
 \end{equation}
 Now we must substitute inequalities (\ref{4w},\ref{3w},\ref{2w}) into (\ref{1w}) and obtain
\begin{equation}
 \mu \left( J_\nu^{(n)}\bigcap B_M \right)
\ge \left( 1-10\sum_{j = m+1}^M \delta_j \right)
 \mu \left( J_\nu^{(n)} \bigcap B_m \right)
 \label{5w}
 \end{equation}
But
$$
\sum_{j = m+1}^M \delta_j = \frac{1}{c} \sum_{j = m+1}^M \frac{1}{j\log j} \le  \frac{1}{c} \log \frac{\log M}{\log m} \le  \frac{1}{c} \log
\left(2+\frac{1}{\gamma}\right)
$$
for $n$ large enough. By (\ref{PARA}) we have
\begin{equation}
1-10\sum_{j = m+1}^M \delta_j \ge 1 -\frac{10}{c}\log \left(2+\frac{1}{\gamma}\right) = \frac{5}{6}.
 \label{6w}
 \end{equation}
 Lemma 2 follows from (\ref{5w},\ref{6w}).

Now we consider a segment $J_\nu^{(n)}$ satisfying (\ref{ini}) and the intersection
\begin{equation}
J_\nu^{(n)} \bigcap B_m = \bigcup_{\kappa=1}^r J_\kappa^{(m)} \neq \varnothing.\label{V}
\end{equation}
By Lemma 2 we know that
$$
\mu \left(\left(\bigcup_{\kappa =1}^r J_\kappa^{(m)} \right)\bigcap B_M \right) \ge \frac{5}{6} \,\, \mu \left( \bigcup_{\kappa =1}^r
J_\kappa^{(m)} \right) .$$

{\bf Lemma 3} \,\,{\it In formula (\ref{V}) there are at least $\lfloor 2r/3\rfloor$ indices $\kappa$ such that
\begin{equation}
\mu\left( J_\kappa^{(m)} \bigcap B_M \right) \ge \frac{1}{2} \mu\left( J_\kappa^{(m)}\right) .\label{VV}
\end{equation}
}

Proof. Define $\kappa$ to be good if (\ref{VV}) holds and define $\kappa$ to be bad otherwise. Let $g$ be the number of good indices and $b$ be
the number of bad indices. Let $\mu =  \mu\left( J_\kappa^{(m)}\right)$. Then by Lemma 2
$$
\frac{5}{6} r \mu\le  \sum_{\kappa =1}^r \mu \left(\left(\bigcup_{\kappa =1}^r J_\kappa^{(m)} \right)\bigcap B_M \right)
=\sum_{\kappa\,\,\text{good}}+ \sum_{\kappa\,\,\text{bad}}.
$$
But
$$
\sum_{\kappa\,\,\text{good}} \le g\mu
$$
and
$$
\sum_{\kappa\,\,\text{bad}} \le \frac{1}{2}b\mu.
$$
Thus
$$
\frac{5}{6} r\mu \le \frac{1}{2} b\mu + g\mu,\,\,\, g+b = r
$$
and so $ g \ge \frac{2r}{3}$. Lemma 3 follows.

{\bf 4.\,\, Proof of Theorem 1. }

We need the following well known result.

 {\bf Theorem (Eggleston \cite{Eg}).}\,\,\ {\it Suppose for every $k$ we have a set
$A_k{=}\bigsqcup\limits_{i=1}^{R_k}I_{k}(i)$ where  $I_{k}(i)$ are segments of the real line of length $|I_{k}(i)|= \Delta_k$. Suppose  that
each interval $I_{k}(i)$ has exactly $N_{k+1}{>}1$ pairwise disjoint subintervals $I_{k+1}(i')$  of length $\Delta_{k+1}$ from the set
 $A_{k+1}$. Let
$R_{k+1}{=}R_k{\cdot} N_{k+1}$ . Suppose that
 $0{<}\nu_0{\le}1$ and  that for every $0{<}\nu{<}\nu_0$  the series
$ \sum_{k=2}^{\infty}\frac{\Delta_{k-1}}{\Delta_k}(R_k(\Delta_k)^\nu)^{-1} $ converges. Then the set $A{=}\bigcap_{k=1}^{\infty}A_k$ has
Hausdorff dimension ${\rm HD}(A){\ge}\nu_0$. }

Suppose we have a segment $J_\nu^{(n)}$ from the union (\ref{union}) with the property (\ref{ini}). Then we have (\ref{V}) and at least $N
=\lfloor 2 r/3\rfloor$ subsegments $J_\kappa^{(m)}$ from (\ref{V}) satisfying  (\ref{VV}). Obviously
$$
N =\lfloor 2 r/3\rfloor = \left\lfloor
 \frac{1}{3} \, \frac{\mu \left( J_\nu^{(n)} \bigcap B_m \right)}{\mu \left( J_\nu^{(m)} \right)}
\right\rfloor
  \ge
  \left\lfloor\frac{1}{3} \,
\frac{\mu \left( J_\nu^{(n)}  \right)}{\mu \left( J_\nu^{(m)} \right)}\right\rfloor =  \left\lfloor\frac{1}{3} \, \frac{2^{l_m}  }{ 2^{l_n}
}\right\rfloor.
$$
Now we can define inductively the sets satisfying Eggleston's theorem. The base of the inductive process is obvious. Assume that for natural $k$
the set $A_k$ consists of segments of the form $ J_\nu^{(n_k)}$ for which (\ref{ini}) holds with $ n = n_k$ ($n_k$ comes from (\ref{NJ})). Then
in each $ J_\nu^{(n_k)}$ one can find exactly
$$
 N_{k+1}  =\left\lfloor
 \frac{1}{3} \, \frac{2^{l_{n_{k+1}}}  }{ 2^{l_{n_k}} }\right\rfloor
 $$
 intervals of the form  $ J_\nu^{(n_{k+1})}$ satisfying
(\ref{ini})  with $ n = n_{k+1}$. Since the series (\ref{ser}) converges  for  $ v \le \frac{\gamma}{\gamma +1}$ (as it was shown in Section 2)
and $l_j$ are defined on (\ref{L}) we see that the series in Eggleston's theorem also converges for the same values of $v$. Hence we deduce that
the Hausdorff dimension of the set $B$ is not less than  $  \frac{\gamma}{\gamma +1}$
 and Theorem 1 is proved.

\vskip+2.0cm

author: Nikolay Moshchevitin

\vskip+0.5cm

e-mail: moshchevitin@mech.math.msu.su, moshchevitin@rambler.ru

\end{document}